\documentclass[11pt,reqno]{amsart}
\usepackage{graphicx}
\usepackage{amssymb}
\usepackage{amsthm}
\usepackage{bm}
\usepackage{amsmath}
\usepackage{mathrsfs}
\usepackage{yfonts}
\usepackage{enumerate}
\usepackage{lscape}
\usepackage{setspace}
\usepackage{hyperref}

\theoremstyle{plain}
    \newtheorem{theorem}{Theorem}
    
    \newtheorem{lemma}[theorem]{Lemma}
    
    \newtheorem{corollary}[theorem]{Corollary}

\theoremstyle{definition} 

\newcommand{\re}{\ensuremath{\mathrm{Re}}}
\newcommand{\im}{\ensuremath{\mathrm{Im}}}

\def\beqs{\begin{eqnarray*}}
\def\eeqs{\end{eqnarray*}}

\def\be{\begin{equation}}
\def\ee{\end{equation}}

\begin{document}
\title[An Improved Upper Bound for the Right-Side Tail of the Crossover Distribution at the Edge of the Rarefaction Fan]{An Improved Upper Bound for the Right-Side Tail of the Crossover Distribution at the Edge of the Rarefaction Fan}

\author[W. Stanton]{William Stanton}
\address{W. Stanton\\
  University of Colorado at Boulder\\
  Boulder, CO\\
  \\
  }
\email{william.stanton@colorado.edu}


\maketitle

This note is a refinement of a calculation done in Sections 5.3 and 5.4 of \cite{CQ}  by Corwin and Quastel.  By improving some of the estimates, we were able to obtain the following result:

\begin{corollary}
Let $T_0 > 0$.  Then there exist $c_1,c_2,c_3<\infty$ depending only on $T_0$ such that for all $T\geq T_0$,
\begin{equation}\label{ldbd}
1 - F_{T,0}^{edge}(s) \leq c_1(e^{-c_2 T^{1/3} s+ e^{-c_3 s^{3/2}}}).
\end{equation}
\end{corollary}

We will need the following useful inequality related to Stirling's approximation.  For $x > 0$,

\begin{equation} \label{Stirling's Inequality}
1<(2\pi)^{{-1/2}}x^{{(1/2)-x}}e^{x}\mathop{\Gamma\/}\nolimits\!\left(x\right)<e^{{1/(12x)}}.
\end{equation}
 
We will also use the fact from Lemma 49 of \cite{CQ} that there exists a constant $C > 0$ such that for all $\re(z) >0$,

\begin{equation}\label{Lemma 49}
\left|{1}/{\Gamma(z)}\right| \leq C e^{2 |z|}.
\end{equation}

The essential improvement of this note over \cite{CQ} is the following lemma, which gives improved bounds for the Airy Upper and Airy Lower Gamma functions.

\begin{lemma}\label{Airy Bounds}
Fix a constant $T_0 > 0$, and let $\kappa_T =  2^{-1/3} T^{1/3}$.  Then there exists a constant $C > 0$ depending only on $T_0$ such that the following inequalities hold:
\begin{enumerate}
\item \begin{equation}\label{Airy Upper Gamma}
|{\rm Ai}^{\Gamma}(x,\kappa_T^{-1},0)| \leq C T^{1/3} \text{ for all $x \in \mathbb{R}$} 
\end{equation}
\item \begin{equation}\label{Airy Lower Gamma Positive}
|{\rm Ai}_{\Gamma}(x,\kappa_T^{-1},0)| \leq C T^{-1/3} e^{-\tfrac{2}{3}x^{3/2}} \text{ for all $x \geq 0$}
\end{equation}
\item \begin{equation}\label{Airy Lower Gamma Negative}
|{\rm Ai}_{\Gamma}(x,\kappa_T^{-1},0)|  \leq C T^{-1/3} e^{2 \kappa_T^{-1} |x|^{1/2}}  \text{ for all $x < 0$}
\end{equation}
\end{enumerate}
\end{lemma}

\emph{Proof of Lemma 2:}

In this proof, $C > 0$ is a constant that can change from line to line, but only depends on $T_0$.
\subsection*{$\mathbf{Ai^{\Gamma}}$ Bounds, $\mathbf{x \geq 0}$}
We begin by proving (\ref{Airy Upper Gamma}).   We begin by proving the $x \geq 0$ case.  

We deform the contour $\tilde{\Gamma}_{\zeta}$ to the vertical line $s_0 + it$, where $s_0 = -\kappa_{T_0}/2$, $t \in (-\infty,\infty)$.    On this new contour, 
\begin{align*}
|\Gamma(\kappa_T^{-1}(s_0 + it))| &= \left|\frac{\Gamma(\kappa_T^{-1}(s_0 + it) + 1)}{\kappa_T^{-1}(s_0 + it)}\right| \\
                                                         &\leq \left|\frac{\Gamma(\re(\kappa_T^{-1}(s_0 + it) + 1))}{\kappa_T^{-1}(s_0 + it)}\right| \\
                                                          &= \left|\frac{\Gamma(\kappa_T^{-1} s_0 + 1)}{\kappa_T^{-1} (s_0 + it)}\right|.
\end{align*}
In the first equation, we use the functional equation $\Gamma(z) = \Gamma(z+1)/z$, and in the second line, we use the simple fact that $|\Gamma(z)| \leq \Gamma(\re(z))$.  
Now, $3/2 \geq \kappa_T^{-1} s_0 + 1 \geq 1/2$, by choice of $s_0$.  Therefore, we conclude that the following is true:
\begin{eqnarray}
|\Gamma(\kappa_T^{-1}(s_0 + it))| &\leq& \frac{\Gamma(3/2)}{|\kappa_T^{-1}(s_0 + i t)|} \\
                                                           &=& \frac{\Gamma(3/2) 2^{-1/3} T^{1/3}}{\sqrt{s_0^2 + t^2}} \\
                                                           &\leq& C T^{1/3}.
\end{eqnarray}

By deforming the contour, we pick up only the residue at $z = 0$, since the poles of $\Gamma(\kappa_T^{-1} z)$ occur at $\kappa_T^{-1} z = -n$, $n \in \mathbb{N} \cup \{0\}$, and
\begin{equation*}
s_0 = - \frac{\kappa_{T_0}}{2} \geq -\frac{\kappa_T}{2} > \kappa_T.
\end{equation*}
Now, it is easy to verify that 
\begin{equation*}
\text{Res}(e^{-1/3 z^3 + xz} \Gamma(\kappa_T^{-1} z),z = 0) = \frac{1}{\kappa_T^{-1}} = 2^{-1/3} T^{1/3}.
\end{equation*}
Therefore, the following is true:
\begin{align*}
|Ai^{\Gamma}(x,\kappa_T^{-1},0)| &\leq 2^{-1/3} T^{1/3} + \left|\int_{-\infty}^{\infty} e^{-\frac{1}{3} (s_0 + it)^3 + x (s_0 + it)} \Gamma(\kappa_T^{-1} (s_0  + it)) dt\right| \\
                                                           &\leq 2^{-1/3} T^{1/3} + C T^{1/3} \int_{-\infty}^{\infty} \left|e^{-\frac{1}{3}(s_0^3 + 3 s_0^2 it - 3s_0 t^2 - i t^3) + x(s_0 + it)}\right| dt \\
                                                           &= 2^{-1/3} T^{1/3} + C T^{1/3} \int_{-\infty}^{\infty} e^{-\frac{-1}{3} s_0^{3} + s_0 t^2 + x s_0} dt \\
                                                           &\leq 2^{-1/3} T^{1/3} + C T^{1/3} \int_{-\infty}^{\infty} e^{-|s_0| t^2} dt \\
                                                           &= C T^{1/3}.
\end{align*}
In the second-to-last line, we use the fact that $x  > 0$ and $s_0 < 0$.

\subsection*{$\mathbf{Ai^{\Gamma}}$ Bound, $\mathbf{x < 0}$}
Now, we move on to the $x < 0$ case.
Let $\tilde{x} = - x$.  We break this proof into three cases.
\subsubsection*{Case 1: $0 < \tilde{x}^{1/2} \leq \kappa_{T_0}$}
We deform the contour $\tilde{\Gamma_{\zeta}}$ to the vertical line $s_0  + i t$, where $s_0 = -\frac{\kappa_{T_0}}{2}$, $t \in (-\infty,\infty)$.  Just as in the $x\geq 0$ case, $|\Gamma(\kappa_T^{-1} (s_0 + it))| \leq C T^{1/3}$.  Furthermore, just as in the $x \geq 0$ case, the only residue picked up when deforming the contour is at $z = 0$.  Therefore,
\begin{align*}
|Ai^{\Gamma}(x,\kappa_T^{-1},0)| &\leq 2^{-1/3} T^{1/3}  + C T^{1/3} \int_{-\infty}^{\infty} e^{\re(-\frac{1}{3}(s_0 + it)^3 + x (s_0 + it))} dt  \\ 
	                                                 &= 2^{-1/3} T^{1/3} + C T^{1/3} \int_{-\infty}^{\infty} e^{-\frac{1}{3} s_0^3 + s_0 t^2 + x s_0} dt \\
					      &\leq 2^{-1/3} T^{1/3} C T^{1/3} \int_{-\infty}^{\infty} e^{\frac{1}{24} \kappa_{T_0}^3 -|\frac{\kappa_{T_0}}{2}| t^2 + \frac{\kappa_{T_0}^3}{2}} \\
                                                            &\leq C T^{1/3}.
\end{align*}
In the second-to-last line, we used the assumption that $\tilde{x} < \kappa_{T_0}^2$.

\subsubsection*{Case 2: $\kappa_{T_0} < \tilde{x}^{1/2} < \kappa_T + 1$}
In the contour integral formula for $Ai^{\Gamma}$, make the change of variables $z = s \tilde{x}^{1/2}$, to find that
\begin{equation*}
Ai^{\Gamma}(x,\kappa_T^{-1},0) = \int_{\tilde{\Gamma_{\zeta}}'} e^{-\tilde{x}^(3/2) ( \frac{1}{3} s^3 + s)} \Gamma(\kappa_T^{-1} \tilde{x}^{1/2} s) \tilde{x}^{1/2} ds,
\end{equation*}
where $\tilde{\Gamma_{\zeta}}'$ is the contour obtained by dividing each point on the contour $\tilde{\Gamma_{\zeta}}$ by $\tilde{x}^{1/2}$.

Deform the contour $\tilde{\Gamma_{\zeta}}'$ to the following contour: a straight line passing from $-\infty e^{-i \frac{3\pi}{4}}$ to $-i$, a semicircular arc passing from $-i$ to $-i$, and a straight line passing from $i$ to $\infty e^{i \frac{3\pi}{4}}$.  Notice that the contour does not pass over any singularities when being deformed.
We deal with the arc first.

Parameterize the arc by $s = e^{i \theta}$, $\theta \in [-\pi/2,\pi/2]$.  First, we bound $\Gamma(\kappa_T^{-1} \tilde{x}^{1/2} s)$ on this arc.  By (\ref{Stirling's Inequality}),
\begin{align*}
|\Gamma(\kappa_T^{-1} \tilde{x}^{1/2} e^{i \theta}+1)| &\leq \exp((\tilde{x}^{1/2} \kappa_T^{-1} \cos(\theta) + 1) \log(\tilde{x}^{1/2} \kappa_T^{-1} \cos(\theta) + 1)) \\
                                                                                          &\leq \exp(((\kappa_T + 1)\kappa_T^{-1} \cos(\theta) + 1) \log((\kappa_T + 1)\kappa_T^{-1} \cos(\theta) + 1) \\
								     &\leq \exp(((1 + \kappa_T^{-1}) \cos(\theta) + 1) \log((1+\kappa_T^{-1}) \cos(\theta) + 1) \\
                 							    &\leq \exp(((1+ \kappa_{T_0}^{-1}) + 1) \log((1+ \kappa_{T_0}^{-1}) + 1) \\
								    &= C.
\end{align*}
In the last inequality, we used the fact that $\kappa_T \geq \kappa_{T_0}$ and $1 \geq \cos(\theta) \geq 0$ for all $\theta \in [-\pi/2,\pi/2]$.

By the functional equation for the gamma function,
\begin{equation*}
|\Gamma(\kappa_T^{-1} \tilde{x}^{1/2} e^{i\theta})| = \frac{|\Gamma(\kappa_T^{-1}\tilde{x}^{1/2} e^{i\theta} + 1)|}{|\kappa_T^{-1} \tilde{x}^{1/2} e^{i\theta}|} \leq \frac{C T^{1/3}}{\tilde{x}^{1/2}},
\end{equation*}
where $C'$ only depends on $\kappa_{T_0}$.

Now, we bound the exponential part of the integrand:
\begin{equation*}
\left|\exp\left(-\tilde{x}^{3/2}\left(\frac{1}{3} (e^{i\theta})^3 + (e^{i\theta})\right)\right)\right| = \exp\left(-\tilde{x}^{3/2} \left(\frac{1}{3} \cos(3 \theta) + \cos(\theta)\right)\right)\leq 1,
\end{equation*}
since $\frac{1}{3} \cos(3\theta) + \cos(\theta) \leq 0$ for all $\theta \in [-\pi/2,\pi/2]$, as is easily verified by basic calculus.

The integrand is therefore bounded by $\frac{C T^{1/3}}{\tilde{x}^{1/2}} \tilde{x}^{1/2} = C T^{1/3}$, and since the arc has length $\pi$, the integral along the arc is bounded by $C T^{1/3}$. 

Now, we move on to the rays.  By symmetry, it is enough to prove the $C T^{1/3}$ bound on the upper ray, parameterized by $s = i + r e^{i \frac{3 \pi}{4}}$, $r \in [0, \infty)$.  We split the argument into two parts.  

The first part is when $r \leq \frac{2}{\sqrt{2}} \left(\frac{1}{\kappa_T^{-1} \tilde{x}^{1/2}} - 1\right)$. Let $z = \kappa_T^{-1} \tilde{x}^{1/2} s$.  Then by choice of $r$, the following is true:
\begin{equation*}
0 \geq Re(z) = \kappa_T^{-1} \tilde{x}^{1/2} (r \cos(\frac{3\pi}{4})) \geq -1 + \kappa_T^{-1} \tilde{x}^{1/2}.
\end{equation*}
Thus, $1 \geq Re(z+1) \geq \kappa_T^{-1} \tilde{x}^{1/2}$.  Therefore,
\begin{align*}
|\Gamma(z)| &= \frac{|\Gamma(z+1)|}{|z|} \\
                     &\leq \frac{2^{-1/3} T^{1/3} |\Gamma(\re(z+1))|}{\tilde{x}^{1/2} |i + r e^{i \frac{3\pi}{4}}|} \\
                     &\leq \frac{2^{-1/3} T^{1/3} \Gamma(\kappa_T^{-1} \tilde{x}^{1/2})}{\tilde{x}^{1/2}} \\
                     &\leq \frac{2^{-1/3} T^{1/3}}{\tilde{x}^{1/2}} \exp(\kappa_T^{-1} \tilde{x}^{1/2} \log (\kappa_T^{-1} \tilde{x}^{1/2}) \\
                     &\leq \frac{2^{-1/3} T^{1/3}}{\tilde{x}^{1/2}} \exp(\kappa_{T}^{-1} (\kappa_T + 1) \log(\kappa_{T}^{-1} (\kappa_T + 1))) \\ 
                    &\leq \frac{C T^{1/3}}{\tilde{x}^{1/2}}.
\end{align*}

The second part is when $r \geq \frac{2}{\sqrt{2}}(\frac{1}{\kappa_T^{-1} \tilde{x}^{1/2}} - 1)$.  Choose a natural number $k$ such that $1/2 \leq \re(z + k) \leq 3/2$.  Note that $k \geq 2$, since by choice of $r$, $\re(z) \leq -1$.  Also, note that for all $j \geq 1$, $|z + j| \geq |\im(z)| \geq 1$, by choice of $r$.  Therefore,  by repeatedly applying the functional equation of the gamma function, we see that the following is true:
\begin{equation*}
|\Gamma(z)| = \frac{|\Gamma(z+k)|}{|z| |z + 1| \cdots |z + k - 1|} \leq \frac{|\Gamma(\re(z + k))|}{|z|} \leq \frac{2^{-1/3} T^{1/3} \Gamma(3/2)}{\tilde{x}^{1/2}|i + r e^{i\frac{3\pi}{4}}|} \leq \frac{C T^{1/3}}{\tilde{x}^{1/2}}.
\end{equation*}

Now, we address the exponential part of the integrand:
\begin{equation*}
|\exp(-\tilde{x}^{3/2} (1/3 s^3 + s))| = \exp(\re(-\tilde{x}^{3/2}(1/3 s^3 + s)))  = \exp\left(-\tilde{x}^{3/2} r^2 -  \tilde{x}^{3/2} r^3 \frac{\sqrt{2}}{6}\right).
\end{equation*}

Therefore, the integral on the ray is bounded above by the following:
\begin{equation*}
\int_{0}^{\infty} \exp(-\tilde{x}^{3/2} r^2) \frac{C T^{1/3}}{\tilde{x}^{1/2}} \tilde{x}^{1/2} dr \leq \frac{C T^{1/3}}{\tilde{x}^{3/4}}.
\end{equation*}

\subsubsection*{Case 3:  $\kappa_T + 1 < \tilde{x}^{1/2}$}
In this argument, we use the same contour as in Case 2.  We prove the $C T^{1/3}$ bound on the semicircular arc first.  We start with bounding the gamma function.  By the functional equation for the gamma function, and (\ref{Stirling's Inequality}),
\begin{equation*}
|\Gamma(\kappa_T^{-1} \tilde{x}^{1/2} e^{i\theta})| \leq  \frac{\exp((\kappa_T^{-1} \tilde{x}^{1/2} \cos(\theta) + 1) \log (\kappa_T^{-1} \tilde{x}^{1/2} \cos(\theta) + 1))}{\kappa_T^{-1} \tilde{x}^{1/2}}.
\end{equation*}

As for the exponential part, it is easy to verify using calculus that as a function of $\theta$, 
\begin{equation*}
\re(-\tilde{x}^{3/2}(\frac{1}{3} s^3 + s) = 
-\tilde{x}^{3/2}(1/3 \cos(3\theta) + \cos(\theta)) 
\end{equation*}
 decreases on $[-\pi/2,0]$ and increases on $[0,\pi/2]$.  The same is clearly true for $(\kappa_T^{-1} \tilde{x}^{1/2} \cos(\theta) + 1) \log (\kappa_T^{-1} \tilde{x}^{1/2} \cos(\theta) + 1)$.  Therefore, the following is true for $s$ on the arc:
\begin{align*}
\left|e^{-x^{3/2}(\frac{1}{3} s^3 + s)} \Gamma(\kappa_T^{-1} \tilde{x}^{3/2} s) \tilde{x}^{1/2}\right| &\leq e^{-\tilde{x}^{3/2}(\frac{1}{3}\cos(\frac{3\pi}{2}) + \cos(\frac{3\pi}{2}))}  \frac{e^{(\kappa_T^{-1} \tilde{x}^{1/2} \cos(\frac{\pi}{2}) + 1) \log (\kappa_T^{-1} \tilde{x}^{1/2} \cos(\frac{\pi}{2}) + 1)}}{\kappa_T^{-1} \tilde{x}^{1/2}} \tilde{x}^{1/2}\\
 &= \frac{1}{\kappa_T^{-1}} \\
&= 2^{-1/3} T^{1/3}.
\end{align*}
Since the arc has length $\pi$, we conclude that the integral along the arc is bounded above by $C T^{1/3}$.

Now, we check the $C T^{1/3}$ bound on the rays.  By symmetry, it is necessary only to check the bound on the upper ray, $s = i + r e^{i \frac{3 \pi}{4}}$.  Choose a natural number $k$ such that $1 \leq Re(z+k) \leq 2$.  Note that $|z+j| \geq Im(z+j) \geq \kappa_T^{-1} \tilde{x}^{1/2} >1$ for all $j \geq 0$.  By again repeatedly applying the functional equation for the gamma function, the following holds:
\begin{equation*}
|\Gamma(\kappa_T^{-1} \tilde{x}^{1/2} (i + r e^{i \frac{3\pi}{4}}))| \leq \frac{|\Gamma(z+k)|}{|z| |z + 1| \cdots |z+k-1|} \leq \frac{\Gamma(2)}{|\kappa_T^{-1}\tilde{x}^{1/2}(i + r e^{i \frac{3 \pi}{4}})|} \leq \frac{C T^{1/3}}{\tilde{x}^{1/2}}.
\end{equation*}

The exponential part is easily bounded, as follows:
\begin{equation*}
\exp(-\tilde{x}^{3/2}(1/3(i + r e^{i \frac{3 \pi}{4}})^3 + (i + r e^{i \frac{3\pi}{4}}))) \leq \exp(- r^2 \tilde{x^{3/2}} - r^3 \tilde{x}^{3/2} \frac{\sqrt{2}}{{6}}).
\end{equation*}
Therefore, the integral along the ray is bounded by the following:
\begin{equation*}
\int_0^{\infty} \exp(-r^2 \tilde{x}^{3/2}) \frac{C T^{1/3}}{\tilde{x}^{1/2}} \tilde{x}^{1/2} dr \leq \frac{C' T^{1/3}}{\tilde{x}^{3/4}}.
\end{equation*}

\subsection*{$\mathbf{{Ai}_{\Gamma}}$ Bound, $\mathbf{x \geq 0}$}

\subsubsection*{Case 1: $0 \leq x^{1/2} \leq \kappa_{T_0}$}  In this case, we deform the contour $\tilde{\Gamma}_{\eta}$ to the contour $z= 1 + r i$, $r \in (-\infty,\infty)$.  The function $1/\Gamma(\kappa_T^{-1} z)$ has no singularities, so we do not pick up any residues.  By (\ref{Lemma 49}) and the functional equation for the gamma function, it now holds that
\begin{equation*}
\left|\frac{1}{\Gamma(\kappa_T^{-1}(1+ ri))}\right|  = \frac{|\kappa_T^{-1} (1 + ri)|}{|\Gamma(\kappa_T^{-1} (1 + ri) + 1)|} \leq  C T^{-1/3} (1 + r) e^{2\kappa_T^{-1} (2 + r^{1/2})}.
\end{equation*}
Therefore, the integral is bounded as follows, with suitable $C, C'$:
\begin{align*}
\int_{-\infty}^{\infty} \exp(1/3(1+ ri)^3 - x (1+ ri)) C T^{-1/3} (1 + r) dr &\leq C T^{-1/3}\int_{-\infty}^{\infty} \exp(1/3 - r^2 - x) (1 + r) e^{2\kappa_T^{-1} (2 + r^{1/2})} dr \\
&\leq C' T^{-1/3} e^{-2/3 x^{3/2}}.
\end{align*}

\subsubsection*{Case 2: $\kappa_{T_0} < x^{1/2}$}  In this case, we make the usual change of variables $z = s x^{1/2}$.  We are free to deform the contour $\tilde{\Gamma}_{\eta}'$ to a steepest descent contour that passes along a straight line from $\infty e^{-i \frac{\pi}{3}}$ to $1$, and then on a straight line from $1$ to $\infty e^{i \frac{\pi}{3}}$.  By Stirling's approximation, it is easy to see that $1/\Gamma(z)$ is bounded above by an absolute constant on both rays.  Therefore, again by the functional equation for the gamma function,
\begin{equation*}
\left|\frac{1}{\Gamma(\kappa_T^{-1} x^{1/2} s)}\right| \leq C T^{-1/3} (1 + r) x^{1/2}.
\end{equation*}

The integral along the upper ray is bounded by the following:
\begin{align*}
\int_0^{\infty} e^{x^{3/2} (\frac{1}{3}(1 + r e^{i\frac{\pi}{3}})^3 + (1 + r e^{i \frac{\pi}{3}}))} C T^{-1/3} (1 + r) x dr &= C T^{-1/3} e^{-\frac{2}{3} x^{3/2}} \int_0^{\infty} e^{-\frac{1}{2} r^2 x^{3/2} - x^{3/2} r^3 + \log(x(1+ r))} dr \\
                  &\leq  C T^{-1/3} e^{-\frac{2}{3} x^{3/2}} \int_0^{\infty} e^{-\frac{1}{2} r^2 x^{3/2}} dr \\
                  &\leq \frac{C T^{-1/3} e^{-\frac{2}{3}x^{3/2}}}{x^{3/4}} 
\end{align*}

\subsection*{$\mathbf{Ai_{\Gamma}}$ Bound, $\mathbf{x < 0}$}

Again, let $\tilde{x} = -x$.   

\subsubsection*{Case 1: $0 \leq \tilde{x} \leq \kappa_{T_0}$}.  We deform the contour $\tilde{\Gamma_{\eta}}$ to the vertical line $1 + r i$, $r \in (-\infty, \infty)$.  We have already shown that 
\begin{equation*}
|\frac{1}{\Gamma(\kappa_T^{-1} z)}| \leq C T^{-1/3} (1 + |r|) e^{4 \kappa_T^{-1} (1 + |r|^{1/2})}
\end{equation*}
 on this contour.  Therefore, the integral is bounded above by the following:
\begin{align*}
\int_{-\infty}^{\infty} e^{Re(1/3 (1 + ri)^3 - x (1 + ri))} C T^{-1/3} (1 + r) e^{4 \kappa_T^{-1} (1 + |r|^{1/2}}dr &= C T^{-1/3} e^{1/3-x} \int_{-\infty}^{\infty} e^{ -r^2} (1 + |r|) e^{4 \kappa_T^{-1} (1 + |r|^{1/2})} dr\\
                                                                                                                                  &\leq C T^{-1/3} e^{1/3 - \kappa_{T_0}} \int_{-\infty}^{\infty} e^{-r^2} (1+ |r|) e^{4 \kappa_T^{-1} (1 + |r|^{1/2})} dr,											
\end{align*}
which is bounded by $C T^{-1/3}$ because $\tilde{x} < \kappa_{T_0}$.

\subsubsection*{Case 2: $\kappa_{T_0} < \tilde{x}^{1/2}$}.  In this case, we first make the change of variables $z = s \tilde{x}^{1/2}$.  Then, we deform the contour to the contour (in the $s$-plane) made up of a straight line passing from $\infty e^{-i \frac{\pi}{4}}$ to $-i$, then a straight line from $-i$ to $i$, and finally a straight line from $i$ to $\infty e^{i\frac{\pi}{4}}$.  

We deal with the vertical line segment from $-i$ to $i$ first.  We parameterize the vertical line segment by $s = -i + ti$, $t \in (0,2)$, so by \ref{Lemma 49}, 

\begin{align*}
\left|\frac{1}{\Gamma(\kappa_T^{-1} \tilde{x}^{1/2} (-i + it))}\right| &= \frac{\left|\kappa_T^{-1} \tilde{x}^{1/2} (-i + it)\right|}{\left|\Gamma(\kappa_T^{-1} \tilde{x}^{1/2} (-i + it) + 1)\right|}  \\
&\leq C T^{-1/3} \tilde{x}^{1/2} e^{2|\kappa_T^{-1} \tilde{x}^{1/2} (-i + it) + 1|}  \\
&\leq C T^{-1/3} \tilde{x}^{1/2} e^{6 \kappa_T^{-1} \tilde{x}^{1/2} }.
\end{align*}

On the vertical line segment, the exponential part is bounded as follows:
\begin{equation*}
|\exp(\tilde{x}^{3/2}((1/3(-i + ti)^3 + (-i + ti)))| = 1.
\end{equation*}
Since the length of the line segment is 2, we conclude that the integral on the line segment is bounded by $C T^{-1/3} \tilde{x}^{1/2} e^{c \tilde{x}} \tilde{x}^{1/2} = C T^{-1/3} \tilde{x} e^{c \kappa_T^{-1} \tilde{x}^{1/2}}$.

Now, we consider the integral along the ray.  By symmetry, we only consider the upper ray.  We parameterize the upper ray by $s = i + r e^{i\frac{\pi}{4}}$, on which it holds that
\begin{equation*}
\left|\frac{1}{\Gamma(\kappa_T^{-1} \tilde{x}^{1/2}(i + r e^{i \frac{\pi}{4}}))}\right| \leq C T^{-1/3} \tilde{x}^{1/2} (1 + r).
\end{equation*}

The exponential part of the integrand is bounded by

\begin{equation*}
 |\exp(\tilde{x}^{3/2}(1/3(i + r e^{i\frac{\pi}{4}})^3 + (i + r e^{i \frac{\pi}{4}})))| \leq \exp(-\tilde{x}^{3/2} r^2) dr,
\end{equation*}

and hence, the integral along the ray satisfies
\begin{align*}
\int_{0}^{\infty} \exp(\tilde{x}^{3/2}(\frac{1}{3}s^3 +s)) C T^{-1/3} \tilde{x}^{1/2} (1 + r) \tilde{x}^{1/2} dr &\leq C T^{-1/3} \tilde{x} \int_0^{\infty} \exp(-\tilde{x}^{3/2} r^2) dr\\
 &= C T^{-1/3} \tilde{x}^{-1/4}.
\end{align*}

\section{Upper Tail of Crossover Distribution}

We give upper bounds for the upper tail of the following \cite{CQ}:

\begin{equation*}
1- F_{T,0}^{\rm edge} ( s) =- \int_{\mathcal{\tilde C}}e^{-\tilde \mu}\frac{d\tilde\mu}{\tilde\mu}[
 \det(I-\tilde{K}_{T,\tilde\mu} )- \det I] ,
\end{equation*}

where $\tilde{K}_{T, \tilde \mu}$ is evaulated on $L^2(s, \infty)$ and 
\begin{equation*}
\tilde{K}_{T,\tilde\mu}(x,y)=\int_{-\infty}^{\infty} \frac{\tilde\mu dt}{e^{-2^{-1/3}T^{1/3}t} - \tilde\mu}  {\rm Ai}^{\Gamma}(x+t ,\kappa_T^{-1},0) {\rm Ai}_{\Gamma}(x+t ,\kappa_T^{-1},0).
\end{equation*}

In order to prove our upper bounds, we follow \cite{CQ} directly, but use the new bounds for the $Ai^{\Gamma}$ and $Ai_{\Gamma}$ we proved in the last section.  We factor the operator $\tilde{K}_{T, \tilde mu}$ into a product of two Hilbert-Schmidt operators, bound them, and then use the continuity of the Fredholm determinant.  We will be consistent with the notation and structure of \cite{CQ} in order to make the argument easy to follow.

By \cite{BS:book}, 

\begin{equation*}
\det(I - \tilde{K}_{T,\tilde\mu}) = \det(I - A),
\end{equation*}

where $A = U^{-1} \tilde{K}_{T,\tilde\mu} U$, and

\begin{equation*}
U f(x) = (x^4 + 1)^{1/2} f(x).
\end{equation*}

We factor $A = A_1 A_2$, where $A_1: L^2(\mathbb{R}) \rightarrow L^2(s, \infty)$, $A_2: L^2(s, \infty) \rightarrow L^2(\mathbb{R})$, and $A_1$ and $A_2$ have the following kernels:
\begin{equation}
A_1(x,t) = Ai^{\Gamma} (x + t, \kappa_T^{-1},0) (x^4 + 1)^{-1/2} (t^4 + 1)^{-1/2}
\end{equation}
\begin{equation}
A_2(t, y) = \frac{\tilde{\mu}}{e^{-\kappa_T t} - \tilde{\mu}} Ai_{\Gamma} (y + t, \kappa_T^{-1}, 0) (x^4 + 1)^{1/2} (t^4 + 1)^{1/2}.
\end{equation}

We will use the fact that since $A_1$ and $A_2$ are Hilbert Schmidt and $A_1 A_2 = A$,

\begin{equation}\label{trbd}
|\det(I+A)-\det I |\le \|A\|_1 e^{ \|A\|_1 +1} \le \|A_1\|_2 \|A_2\|_2 e^{\|A_1\|_2 \|A_2\|_2 +1}.
\end{equation}

By (\ref{Airy Upper Gamma}) and the definition of the Hilbert-Schmidt norm, the following is true:
\begin{align}
\|A_1\|_2^2 &= \int_{s}^{\infty} \int_{-\infty}^{\infty} (Ai^{\Gamma}(x + t, \kappa_T^{-1},0) (x^4 + 1)^{-1/2} (t^4 + 1)^{-1/2})^2 dx dt  \\
                     &\leq   C T^{2/3} \int_{-\infty}^{\infty} \int_{-\infty}^{\infty}  (x^4 + 1)^{-1/2} (t^4 + 1)^{-1/2})^2 dx dt  \\
                     &= C T^{2/3}.
\end{align}

In order to bound $A_2$, we will use the following bound, along with the bound we proved for $Ai_{\Gamma}$:
\begin{equation}
\left|\frac{\tilde{\mu}}{e^{-\kappa_T t} - \tilde{\mu}}\right| \leq C |\tilde{\mu}|(e^{2\kappa_T t} \wedge 1),
\end{equation}
which is formula (116) in \cite{CQ}.

By our upper bounds for $Ai^{\Gamma}$, the following is true:
\begin{align*}
\|A_2\|_2^2 &\leq C T^{-2/3} |\tilde{\mu}|^2 \left( \int_s^{\infty} dx \int_{-x}^{\infty} dt (e^{2 \kappa_T t} \wedge 1) e^{-\frac{4}{3} |x + t|^{3/2}} (x^4 + 1)(t^4 + 1) \right. \\
                     & \left. + \int_s^{\infty} dx \int_{-\infty}^{-x} dt (e^{2 \kappa_T t} \wedge 1) |x+t| e^{2 \kappa_T^{-1} |x+t|^{1/2}} (x^4 + 1) (t^4 + 1) \right)  C |\tilde{\mu}|^2 (I_1 + I_2).
\end{align*}

First of all, notice that since $s >>1$, $t < 0$ over the range of the inner integral.  Thus, $e^{2 \kappa_T t} \wedge 1 = e^{2 \kappa_T t}$ in the inner integral.  Make the change of variables $u = t/x$ in the inner integral.  Then the integral now has the following form:
\begin{equation*}
\int_{-\infty}^{-1} e^{2 \kappa_T x u}  |x+x u| e^{2 \kappa_T^{-1} x^{1/2} |1 + u|} (x^4 + 1) (t^4 + 1) x du.
\end{equation*}

Now, choose $s > 64 \kappa_{T_0}^{-4}$.  Then since $x \geq s$,  it is easy to verify that 
\begin{equation*}
-2 \kappa_T x + 8 \kappa_T^{-1} |x|^{1/2} < -\kappa_T x
\end{equation*}
This implies the following inequality:
\begin{align*}
I_2 &\leq  \int_s^{\infty} dx \int_{-\infty}^{-1}x  du e^{- \kappa_T x u} |x + ux| (x^4 + 1) ((xu)^4 + 1) \\
        &\leq C e^{-\kappa_T s} s^{9}.
\end{align*}

Now, we bound $I_1$.  We write $I_1 = I_3 + I_4$, where the $t$-integration runs from $-x$ to $-x/2$ in $I_3$, and the $t$-integration runs from $-x/2$ to $\infty$ in $I_4$.  Now,
\begin{equation*}
I_3 \leq \int_s^{\infty} dx \int_{-x}^{-x/2} dt e^{2 \kappa_T t} (x^4 + 1)(t^4 + 1) \leq C s^8 e^{-\kappa_T s},
\end{equation*}
and
\begin{equation*}
I_4 \leq \int_s^{\infty} dx \int_{-x/2}^{\infty} dt e^{-4/3|x+t|^{3/2}} (x^4 + 1)(t^4 + 1) \leq C s^{8} e^{-c s^{3/2}}.
\end{equation*}
Therefore,
\begin{equation*}
\|A_2\|_2^2 \leq C |\tilde{\mu}|^2 s^9 (e^{-\kappa_T s} + e^{-c s^{3/2}}).
\end{equation*}

Now, we combine our estimates for $A_1$ and $A_2$, use \ref{trbd}, and conclude the following:
\begin{equation*}
|\det(I - \tilde{K}_{T, \tilde{\mu}})| \leq C \|\tilde{\mu}\| s^{9/2} (e^{-\kappa_T s} + e^{-c s^{3/2}}) e^{\frac{1}{2}\tilde{\mu}|},
\end{equation*}
where we have chosen $s$ large enough that $C s^{9/2} (e^{-\kappa_T s} + e^{-c s^{3/2}}) \leq 1/2$. 

Now, we insert this estimate into the formula for the probability distribution:
\begin{align*}
|1 - F_{T,0}^{edge}(s)| &\leq  \int_{\tilde{C}} |e^{-\tilde{\mu}}| \frac{d |\tilde{\mu}|}{|\tilde{\mu}|} \|\tilde{\mu}\| s^{9/2} (e^{-\kappa_T s} + e^{-c s^{3/2}}) e^{\frac{1}{2} \tilde{\mu}}  \\
                                       &\leq C s^{9/2} (e^{-c T^{1/3} s} + e^{-c s^{3/2}}).
\end{align*}
In the second inequality, we have used the fact that $|\int_{\tilde{C}} e^{-\tilde{\mu}} |\tilde{\mu}| \frac{d |\tilde{\mu}|}{|\tilde{\mu}|} e^{\frac{1}{2} |\tilde{\mu}|}|$ converges.  

Now, $T \geq T_0$, so for all sufficiently large $s$, there exists $c' > 0$ such that the following is true:
\begin{align*}
s^{9/2} (e^{-c T^{1/3} s} + e^{-c s^{3/2}}) &\leq e^{-c T_0^{1/3} s + \log s^{9/2}} + e^{-c s^{3/2} + \log s^{9/2}} \\
                                                                          &\leq e^{-c' T_0^{1/3} s} + e^{-c' s^{3/2}}.
\end{align*}
Therefore, we conclude the following:
\begin{equation*}
|1 - F_{T,0}^{edge}(s)| \leq  C  (e^{-c' T^{1/3} s} + e^{-c' s^{3/2}}).
\end{equation*}

\bibliographystyle{alpha}

\end{document}